\documentclass[11pt]{article}

\usepackage{amsfonts,amssymb,amsmath,amsthm,epsfig,euscript}
\usepackage{epic, eepic, graphics}
\usepackage{graphicx}

\newcommand{\bob}[1]{m_{#1}}
\newcommand{\rank}{\ell} % rank (but not in the normal poset terminology)
\newcommand{\srank}{\ell^{\star}} % special rank
\newcommand{\myadd}{\Phi}
\newcommand{\remove}{\Psi}
\newcommand{\addi}{\textsf{Add1}}
\newcommand{\addii}{\textsf{Add2}}
\newcommand{\addiii}{\textsf{Add3}}
\newcommand{\subi}{\textsf{Sub1}}
\newcommand{\subii}{\textsf{Sub2}}
\newcommand{\subiii}{\textsf{Sub3}}
\newcommand{\M}{\mathcal{M}}

\newtheorem{theorem}{Theorem}

\newtheorem{lemma}[theorem]{Lemma}

\newtheorem{example}{Example}

\newtheorem{conjecture}{Conjecture}

\long\def\symbolfootnote[#1]#2{\begingroup
\def\thefootnote{\fnsymbol{footnote}}\footnote[#1]{#2}\endgroup}

\newcommand{\tpt}{$(\mathbf{2+2})$}
\newcommand{\tptp}{$\mathbf{2+2}$}
\newcommand{\tpo}{$(\mathbf{3+1})$}

\def\A{\mathcal{A}}

\def\P{\mathcal{P}}

\def\N{\mathbb{N}}

\newcommand{\pattern}{
  \begin{minipage}[c]{1.45em}\scalebox{0.5}{\includegraphics{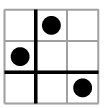}}
  \end{minipage}
}

\DeclareMathOperator{\length}{\mathrm{length}}
\DeclareMathOperator{\asc}{\mathrm{asc}}
\DeclareMathOperator{\minmax}{\mathrm{minmax}}
\DeclareMathOperator{\levels}{\mathrm{levels}}
\DeclareMathOperator{\lds}{\mathrm{lds}}
\DeclareMathOperator{\size}{\mathrm{size}}
\DeclareMathOperator{\run}{\mathrm{run}}
\DeclareMathOperator{\last}{\mathrm{last}}
\DeclareMathOperator{\zeros}{\mathrm{zeros}}

\title{Enumerating \tpt-free posets by the number of minimal elements and other statistics}

\author{
Sergey Kitaev\footnote{The work presented here was supported by grant no. 090038011 from the Icelandic Research Fund.} \\
\small The Mathematics Institute\\[-0.8ex]
\small School of Computer Science\\[-0.8ex]
\small Reykjav\'{i}k University \\[-0.8ex]
\small IS-103 Reykjav\'{i}k, Iceland\\[-0.8ex]
\small \texttt{sergey@ru.is} \and
Jeffrey Remmel\footnote{Partially supported by NSF grant DMS 0654060.} \\
\small Department of Mathematics\\[-0.8ex]
\small University of California, San Diego\\[-0.8ex]
\small La Jolla, CA 92093-0112. USA\\[-0.8ex]
\small \texttt{remmel@math.ucsd.edu}
}

\date{\small Submitted: Date 1;  Accepted: Date 2;
 Published: Date 3.\\
\small MR Subject Classifications: 05A15}

\begin{document}
\maketitle

\begin{abstract}
An unlabeled poset is said to be \tpt-free if it does not contain an
induced subposet that is isomorphic to {\tptp}, the union of two
disjoint 2-element chains. Let $p_n$ denote the number of \tpt-free
posets of size $n$. In a recent paper, Bousquet-M\'elou et
al.~\cite{BCDK} found, using so called ascent sequences, the
generating function for the number of \tpt-free posets of size $n$:
$P(t)=\sum_{n \geq 0} p_n t^n = \sum_{n\geq 0} \prod_{i=1}^{n}
\left( 1-(1-t)^i \right)$. We extend this result in two ways. First,
we find the generating function for \tpt-free posets when four
statistics are taken into account, one of which is the number of
minimal elements in a poset. Second, we show that if $p_{n,k}$
equals the number of \tpt-free posets of size $n$ with $k$ minimal
elements, then $P(t,z)=\sum_{n,k \geq 0} p_{n,k} t^n z^k = 1+
\sum_{n \geq 0} \frac{zt}{(1-zt)^{n+1}} \prod_{i=1}^n (1-(1-t)^i)$.
The second result cannot be derived from the first one by a
substitution. On the other hand, $P(t)$ can easily be obtained from
$P(t,z)$ thus providing an alternative proof for the enumeration
result in~\cite{BCDK}. Moreover, we conjecture a simpler form of
writing $P(t,z)$. Our enumeration results are extended to certain
restricted permutations and to regular linearized chord diagrams
through bijections in~\cite{BCDK,cdk}. Finally, we define a subset
of ascent sequences counted by the Catalan numbers and we discuss
its relations with \tpt- and \tpo-free posets.
\end{abstract}

\section{Introduction}

An unlabeled poset is said to be \emph{{\tpt}-free} if it does not
contain an induced subposet that is isomorphic to {\tptp}, the union
of two disjoint 2-element chains. We let $\P$ (resp. $\P_n$) denote
the set of \tpt-free posets (resp. on $n$ elements).
Fishburn~\cite{fishburn} showed that a poset is {\tpt}-free
precisely when it is isomorphic to an {\em interval order}.
Another important characterization of \tpt-free posets,
set \cite{FISH_BOOK,FISH_OPER, SKANDERA}, is that
a poset is {\tpt}-free if
and only if the collection of strict principal down-sets
can be linearly ordered by inclusion. Here for  any poset $\mathbf{P} =
(P,<_p)$ and $x \in P$, the strict principal
down set of $x$, $D(x)$, in $\mathbf{P}$ is
the set of all $y \in P$ such that $y <_p x$.  The {\em trivial down-set} is the empty set. Thus if $\mathbf{P}$ is a {\tpt}-free poset, we can
write $D(\mathbf{P})= \{D(x):x \in P\}$ as
$$D(\mathbf{P}) =\{D_0,D_1, \ldots ,D_k\}$$
where $\emptyset = D_0 \subset D_1 \subset \cdots \subset D_k$. In
such a situation, we say that $x \in P$ has level $i$ if $D(x) = D_i$.

Let $p_n$ be the number of of \tpt-free posets on $n$ elements.
El-Zahar~\cite{ZAHAR} and Khamis~\cite{smkhamis} used a
recursive description of \tpt-free posets
to derive a pair of functional equations that define
the series $P(t)$. However, they did not solve these equations.
Haxell, McDonald and Thomasson~\cite{haxell} provided an algorithm,
based on a complicated recurrence relation, to produce the first few values
of $p_n$.  Bousquet-M\'elou et al.~\cite{BCDK} showed that the generating
function for the number $p_n$ of \tpt-free posets on $n$ elements is
\begin{equation}\label{gf}
P(t)= \sum_{n\geq 0} p_n \, t^n=\sum_{n\ge 0} \prod_{i=1}^{n} \left(
1-(1-t)^i\right). \end{equation} Note that the term corresponding to
$n=0$ in the last sum is 1.

Zagier~\cite{zagier} proved that (\ref{gf}) is also the generating
function which counts certain involutions introduced by
Stoimenow~\cite{stoim}. Bousquet-M\'elou et al.~\cite{BCDK} gave a
bijections between \tpt-free posets and such involutions,
between \tpt-free posets and a certain restricted class of  permutations,
and between \tpt-free posets and {\em ascent
sequences}. A sequence $(x_1,\dots , x_n ) \in \mathbb{N}^n$ is an
{\em{ascent sequence of length $n$}} if and only if it satisfies $x_1=0$ and
$x_i \in [0,1+\asc(x_1,\dots , x_{i-1})]$ for all $2\leq i \leq n$.
Here, for any integer sequence $(x_1,\dots , x_i)$, the number of
{\em{ascents}} of this sequence is
\begin{equation*}
  \asc(x_1,\dots , x_{i}) = |\{\,1\leq j <i\,:\, x_j<x_{j+1}\,\}|.
\end{equation*}
For instance, (0, 1, 0, 2, 3, 1, 0, 0, 2) is an ascent sequence. We
let $\A$ denote the set of all ascent sequences where we assume the empty
word is also an ascent sequence.

To define the bijection between \tpt-free posets and ascent
sequences, Bousquet-M\'elou et al.~\cite{BCDK} used a step by step
decomposition of a \tpt-free poset $P$ where are each step one
removes a maximal element located on the lowest level together with
certain relations. If one records the levels from which one removed
such maximal elements and then reads the resulting sequence
backwards, one obtains an ascent sequence associated to the poset.
We shall give a detailed account to this bijection in
Section~\ref{decomp}. In the process of decomposing the \tpt-free
poset $P$, one will reach a point where the remaining poset consists
of an antichain, possibly having one element. We define $\lds(P)$ to
be the maximum size of such an antichain, which is also equal to the
size of the down-set of the last removed element that has a
non-trivial down-set. By definition, the value of $\lds$ on an
antichain is 0 as there are no non-trivial down-sets for such a
poset.

Bousquet-M\'elou et al.~\cite{BCDK} studied a
more general generating function $F(t,u,v)$ of \tpt-free posets
according to $\size$=``number of elements'' (variable $t$),
$\levels$=``{\em number of levels}'' (variable $u$),
and $\minmax$=``level of minimum maximal element'' (variable $v$).
The first few terms of $F(t,u,v)$ are
$$
F(t,u,v)=1+t+(1+uv)t^2+ (1+2uv+u+u^2v^2)t^3+ O(t^4).
$$
An explicit form of $F(t,u,v)$ can be obtained from~\cite[Lemma
13]{BCDK} and~\cite[Proposition 14]{BCDK}. The key strategy used by
Bousquet-M\'elou et al.~\cite{BCDK} to derive such formulas was to
translate the appropriate statistics on \tpt-free posets to
statistics on ascent sequences since it is much easier to count
ascent sequences.

The main result of this paper, Theorem~\ref{mainres}, is an explicit
form of the generating function $G(t,u,v,z,x)$ for a generalization
of $F(t,u,v)$, when two more statistics are taken into account ---
$\min$=``number of minimal elements'' in a poset (variable $z$) and
$\lds$=``size of {\em non-trivial last down-set}'' (variable $x$).
That is, we shall find an explicit formula for
$$G(t,u,v,z,x)=\sum_{P\in
\P}t^{\size(P)}u^{\levels(P)}v^{\minmax(P)}z^{\min(P)}x^{\lds(P)}$$
where, as above, $\P$ is the set of all \tpt-free posets. As in
~\cite{BCDK}, to find $G(t,u,v,z,x)$, we translate our problem on
\tpt-free posets to an equivalent problem on ascent sequences. That
is,  we define the following statistics on an ascent sequence:
$\length$=``the number of elements in the sequence,'' $\last$=``the
rightmost element of the sequence,'' $\zeros$=``the number of 0's in
the sequence,'' $\run$=``the number of elements in the leftmost run
of 0's''=``the number of 0's to the left of the leftmost non-zero
element.'' By definition, if there are no non-zero elements in an
ascent sequence, the value of $\run$ is 0. Then we shall prove the
following.

\begin{lemma} The function $G(t,u,v,z,x)$ defined above can
alternatively be defined on ascent sequences as
\begin{equation}
G(t,u,v,z,x)=\sum_{w\in
\A}t^{\length(w)}u^{\asc(w)}v^{\last(w)}z^{\zeros(w)}x^{\run(w)}=\sum_{n,a,\ell,m,r\geq
0}G_{n,a,\ell,m,r}t^nu^av^{\ell}z^mx^r.
\end{equation}
\end{lemma}
\begin{proof}
To prove the statement we need to show equidistribution of the statistics involved. All but one case follow from the results
in~\cite{BCDK}. More precisely, we can use the bijection from
\tpt-free posets to ascent sequences presented in~\cite{BCDK} which
sends $\size\rightarrow\length$, $\levels\rightarrow \asc$,
$\minmax\rightarrow\last$, and $\min\rightarrow\zeros$.
The fact that $\lds$ goes to $\run$ follows from the bijection.
That is, in the process of decomposing the poset, there
will be a point where we remove the element, say $e$, whose down-set
gives $\lds$. At that point, we
will be left with incomparable elements located on
level 0, which gives in the corresponding ascent sequence the
initial run of 0's followed by 1 corresponding to $e$ located on
level 1.
\end{proof}
We shall also give an explicit form of a specialization of
$G(t,u,v,z,x)$, namely $G(t,1,1,z,1)$, which cannot be derived
directly from $G(t,u,v,z,x)$ by the substitution. More precisely,
let $p_{n,k}$ denote the number of \tpt-free posets of size $n$ with
$k$ minimal elements or, equivalently, the number of ascent
sequences of length $n$ with $k$ zeros. Then we shall prove that
\begin{equation}
P(t,z) = \sum_{n,k \geq 0} p_{n,k}t^n z^k =
1+ \sum_{n \geq 0} \frac{zt}{(1-zt)^{n+1}} \prod_{i=1}^n (1-(1-t)^i).
\end{equation}
Moreover, we will conjecture a simpler form of writing $P(t,z)$ (see
Conjecture~\ref{conject}).

A poset $P$ is \tpo-free if it does not contain, as an induced
subposet, a 3-element chain and an element which incomparable to the
elements in the 3-element chain. It is known that the number of
posets avoiding \tpt\ and \tpo\ is given by the Catalan numbers
(see~\cite{stanley,SKANDERA}).
 Define a {\em restricted ascent sequence} as follows. A sequence
$(x_1,\dots , x_n ) \in \mathbb{N}^n$ is a restricted ascent
sequence of length $n$ if it satisfies $x_1=0$ and $x_i \in
[m-1,1+\asc(x_1,\dots , x_{i-1})]$ for all $2\leq i \leq n$, where
$m$ is the maximum element in $(x_1,\dots , x_{i-1})$. For instance,
(0, 1, 0, 2, 3, 2, 2, 3, 2) is a restricted ascent sequence, whereas
(0, 1, 0, 2, 0, 1) is not. Thus, the difference here from the
definition of an ascent sequence is 0 substituted by $m-1$. We shall
show that restricted ascent sequences are counted by the Catalan
numbers. For $n \leq 6$, the bijection in~\cite{BCDK} sends
\tpt- and \tpo-free posets to restricted ascent sequences
which lead us to initially conjecture that it always
the case that the bijection in~\cite{BCDK} sends
\tpt- and \tpo-free posets to restricted ascent sequences. However,
this is not true as we shall produce counter examples when
$n =7$.

This paper is organized as follows. In Section~\ref{decomp}, we
follow \cite[Section 3]{BCDK} to describe a \tpt-free posets
decomposition that gives a bijection between \tpt-free posets and
ascent sequences. The bijection allows us to reduce the enumerative
problem on posets to that on ascent sequences. In
Section~\ref{mainresults} we find explicitly the function
$G(t,u,v,z,x)$ using the ascent sequences (see
Theorem~\ref{mainres}). In Section~\ref{lenzeros}, we shall derive
our formula for $P(t,z)$ and state a conjecture on a different form
for it. We also show in Section~\ref{lenzeros} how to get $P(t)$
from $P(t,z)$ thus providing an alternative proof for the
enumeration in~\cite{BCDK}. Finally, in Section~\ref{catalan} we
define a subset of ascent sequences counted by the Catalan numbers
and discuss its relations to \tpt- and \tpo-free posets.

\section{{\tpt}-free posets and ascent sequences}\label{decomp}

In this section, we shall review the bijection
between \tpt-free posets and ascent sequences
given in \cite[Section 3]{BCDK}.
In order to do this, Bousquet-M\'elou et al.~\cite{BCDK} introduced
 two operations on posets
in $\P_n$. The first is an addition operation; it adds an element to
$P \in \P_n$ that results in $Q \in \P_{n+1}$. The second is a
removal operation; it removes a maximal element $\bob{P}$ from $P
\in \P_{n}$ and results in $Q\in\P_{n-1}$. Before giving these
operations we need to define some terminology.

Let $D(x)$ be the set of predecessors of $x$ (the strict down-set of
$x$): $D(x)=\{\, y : y<x \,\}$. Clearly, any poset is uniquely
specified by listing the sets of predecessors. It is
well-known---see for example Khamis~\cite{smkhamis}---that a poset
is {\tpt}-free if and only if its sets of predecessors, $\{D(x) :
x\in P\}$, can be linearly ordered by inclusion. Let
$$D(P) = (D_0,D_1,\dots,D_{k-1})
$$ with $D_0\subset D_1\subset\dots\subset D_{k-1}$ be that chain. In
this context we define $D_i(P)=D_i$ and $\rank(P) = k$. We say the
element $x$ is at \emph{level} $i$ in $P$ if $D(x)=D_i$ and we write
$\rank(x)=i$ .  The set of all elements at level $i$ we denote
$L_i(P) = \{\, x\in P : \rank(x)=i \,\}$ and we let
$$L(P) = \big(\,L_0(P),L_1(P),\dots,L_{k-1}(P)\,\big).
$$
For instance, $L_0(P)$ is the set of minimal elements and
$L_{k-1}(P)$ is the set of maximal elements whose set of
predecessors is also maximal.
% The empty level $\rank(P)$ is also valid.
Let $\bob{P}$ be a maximal element of $P$ whose set of predecessors
is smallest. This element may not be unique but the level on which
it resides is.
% Since we are dealing with unlabeled structures it does not
% matter which one we select.
Let us write $\srank(P)=\rank(\bob{P})$.

\begin{example}
Consider the {\tpt}-free poset $P$:
$$
\begin{minipage}{5.2em}
  \includegraphics[scale=0.6]{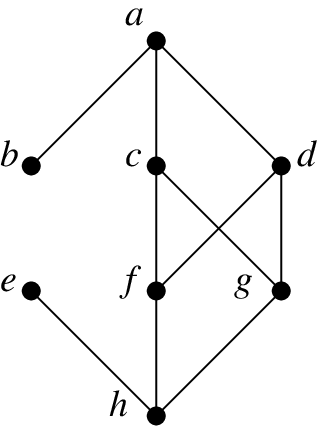}
\end{minipage}
\quad = \quad
\begin{minipage}{4em}
  \includegraphics[scale=0.6]{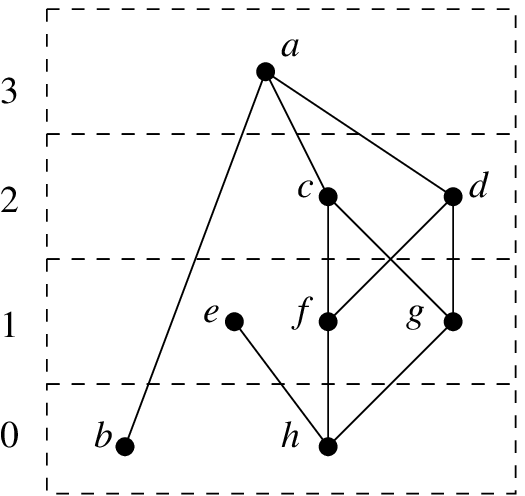}
\end{minipage}
\qquad\quad
$$
The diagram on the right shows the poset redrawn according to the
level numbers of the sets of predecessors. We have $D(a)
=\{b,c,d,f,g,h\}$, $D(b) =\emptyset$, $D(c) = D(d) = \{f,g,h\}$,
$D(e) = D(f) = D(g) =\{h\}$ and $D(h) =\emptyset$. These may be
ordered by inclusion as
$$
\begin{array}{ccccccc}
  \underbrace{D(h) = D(b)} &\!\!\!\subset\!\!\!&
  \underbrace{D(e) = D(f) = D(g)} &\!\!\!\subset\!\!\!&
  \underbrace{D(c) = D(d)} &\!\!\!\subset\!\!\!&
  \underbrace{D(a)}. \\[1.8ex]
  \rank(h)=\rank(b)=0 &&
  \rank(e)=\rank(f)=\rank(g)=1 &&
  \rank(c)=\rank(d)=2 &&
  \rank(a)=3
\end{array}
$$ Thus $\rank(P)=4$. The maximal elements of $P$ are $e$ and
$a$. Since $D(e) \subset D(a)$ we have $\bob{P}=e$ and
$\srank(P)=1$. In addition, $D_0=\emptyset$, $D_1=\{h\}$,
$D_2=\{f,g,h\}$ and $D_3=\{b,c,d,f,g,h\}$. With $L_i = L_i(P)$ we
also have $L_0=\{h,b\}$, $L_1=\{e,f,g\}$, $L_2=\{c,d\}$ and
$L_3=\{a\}$.
%% \begin{align*}
%% L_3(P)&=\{a\}     & D_3(P)&=\{b,c,d,f,g,h\} \\
%% L_2(P)&=\{c,d\}   & D_2(P)&=\{f,g,h\} \\
%% L_1(P)&=\{e,f,g\} & D_1(P)&=\{h\}\\
%% L_0(P)&=\{h,b\}   & D_0(P)&=\emptyset.
%% \end{align*}
\end{example}

Clearly, any {\tpt}-free poset $P$ is determined by the pair
$\big(D(P),L(P)\big)$. Thus when defining the addition and
subtraction operations below it suffices to specify how $D(P)$ and
$L(P)$ change.

%The proof consists of showing how each $P \in \P_{n+1}$ is made by adding a new element to
%a particular $Q \in \P_n$.
The addition operation is actually one of three addition operations,
which will depend on a parameter of $P$. These addition operations
are, in a sense, disjoint.  The first addition operation will result
in $\rank(P)=\rank(Q)$ whereas the second two addition operations
will result in $\rank(Q)=\rank(P)+1$.

Given $P \in \P_n$ and $0\leq i \leq \rank(P)$, let $\myadd(P,i)$ be
the poset $Q$ obtained from $P$ according to the following:
\begin{enumerate}
\item[(\addi)] If $0\leq i \leq \srank(P)$, then introduce a
  new maximal element $z$ on level $i$ which covers the same elements
  as the other elements on level $i$. In terms of predecessors and levels,
  $D(Q)=D(P)$ and
%  $D_j(Q) = D_j(P)$ for all $0\leq j< \rank(P)$ and
%  $D(z)=D_i(P)$.
  \begin{equation*}
    L_j(Q) =
    \begin{cases}
      L_j(P)            & \text{if $j\neq i$}, \\
      L_i(P) \cup \{z\} & \text{if $j=i$}.
    \end{cases}
  \end{equation*}

\item[(\addii)]  If $i = 1 + \rank(P)$, then add a new
element covering all the maximal elements of $P$.

\item[(\addiii)] If $\srank(P)<i\leq\rank(P)$, then let
  $\M$ be the set of maximal elements $x\in P$ with
  $\rank(x)<i$.  Introduce a new element $z$ and set
  $D(z)=D_i(P)$. For all elements of $P$ on level $i$ and above,
  ensure they are greater than every element in $\M$. In
  terms of predecessors and levels,
  \begin{align*}
    D_j(Q) &=
    \begin{cases}
      D_j(P)                       & \text{if $0\leq j \leq i$},\\
      D_{j-1}(P)\cup \M & \text{if $i<j\leq \rank(P)$},
    \end{cases} \\
    \text{and} \\
    L_j(Q) &=
    \begin{cases}
      L_j(P)     & \text{if $0\leq j <i$}, \\
      \{z\}      & \text{if $j=i$},\\
      L_{j-1}(P) & \text{if $i<j\leq \rank(P)$}.
    \end{cases}
  \end{align*}

\end{enumerate}

An important property of the above addition operations is that
$\srank(\myadd(P,i))=i$, since all maximal elements below level $i$
are covered and therefore not maximal in $\myadd(P,i)$.  This allows
us to give the three rules for reversing each of the addition rules
above.

As before, let $\bob{P}$ be a maximal element of $P$ whose set of
predecessors is smallest.  For non-empty $P\in\P_n$, let
$\remove(P)=(Q,i)$ where $i=\srank(P)$ and $Q$ is the poset that
results from applying:
%Deleting an element
\begin{enumerate}
\item[(\subi)] If $\bob{P}$ is not alone on level $i$,
  then remove $\bob{P}$. In terms of predecessors and levels,
  $D(Q) = D(P)$ and
  \begin{equation*}
    L_j(Q) =
    \begin{cases}
      L_j(P)               & \text{if $j\neq i$}, \\
      L_i(P) - \{\bob{P}\} & \text{if $j=i$}.
    \end{cases}
  \end{equation*}
%  This operation is illustrated by reversing the diagram in Figure~\ref{addition_one}.
  \smallskip

\item[(\subii)] If $\bob{P}$ is alone on level
  $i= \rank(P)$, then remove the unique element of level $i$.

\item[(\subiii)] If $\bob{P}$ is alone on level
  $i \leq \rank(P)-1$, then set
  $\N=D_{i+1}(P)-D_{i}(P)$.  Make each element in
  $\N$ a maximal element of the poset by removing any
  covers. Finally, remove the element $\bob{P}$.  In terms of
  predecessors and levels,
  \begin{align*}
    D_j(Q) &=
    \begin{cases}
    D_j(P)                   & \text{if $0\leq j < i$},\\
    D_{j+1}(P) - \N & \text{if $i\leq j < \rank(P)-1$},
    \end{cases} \\
    \text{and} \\
    L_j(Q) &=
    \begin{cases}
    L_j(P)     & \text{if $0\leq j <i$}, \\
    L_{j+1}(P) & \text{if $i\leq j < \rank(P)-1$}.
    \end{cases}
  \end{align*}

\end{enumerate}

We provide an example showing how to find the ascent sequence
corresponding to a given unlabeled poset. See \cite[Section 3]{BCDK}
for an example of how to construct an unlabeled {\tpt}-free poset
corresponding to a given ascent sequence.

\begin{example}\label{example}
Let $P$ be the unlabeled {\tpt}-free poset with this Hasse diagram:
$$
\begin{minipage}{8.6em}
  \includegraphics[scale=0.85]{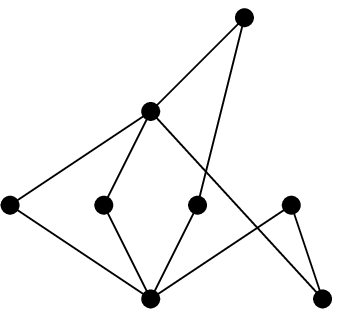}
\end{minipage} =\quad
\begin{minipage}{10em}
  \scalebox{0.70}{\includegraphics{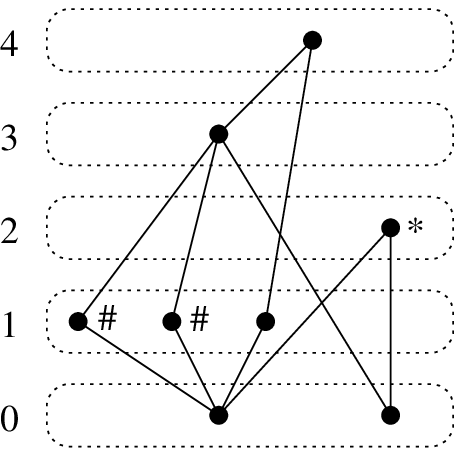}}
\end{minipage}
$$
The diagram on the right shows the poset redrawn according to the
level numbers of the sets of predecessors.  The element $\bob{P}$ is
marked by $*$ and $\rank(\bob{P})=\srank(P)=2$ so $x_8=2$.  Since
$\bob{P}$ is alone on level $2<\rank(P)$ we apply rule {\subii} to
remove it.  The elements corresponding to $\N$ are indicated by
\#'s.  For each of the elements above level 2, destroy any relations
to these \# elements.  Remove $\bob{P}$. Adjust the
level numbers accordingly.\medskip\\
$$\mapsto\quad
\begin{minipage}{10em}\scalebox{0.70}
  {\includegraphics{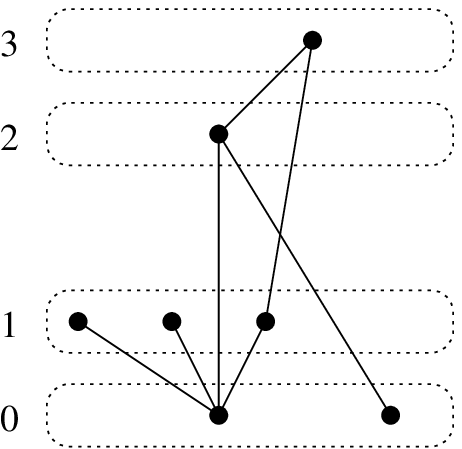}}
\end{minipage} =\quad
\begin{minipage}{10em}
  \scalebox{0.70}{\includegraphics{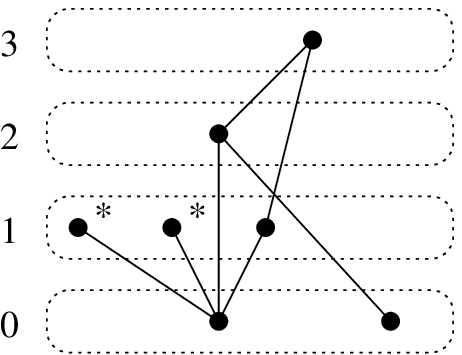}}
\end{minipage}
$$
There are now two elements marked by $*$ that may be considered
$\bob{P}$, however only the level number is important. Thus $x_7=1$
and remove either of the points according to rule {\subi} since
there is more than one element on level 1. Repeat again to give
$x_6=1$.
$$\mapsto\quad
\begin{minipage}{8.5em}
  \scalebox{0.70}{\includegraphics{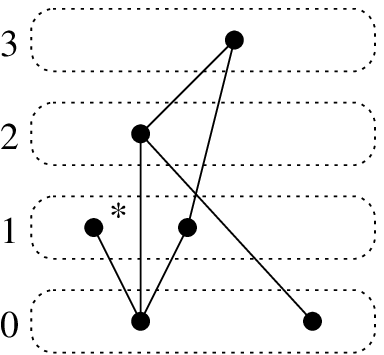}}
\end{minipage}
\mapsto \quad
\begin{minipage}{7.4em}
  \scalebox{0.70}{\includegraphics{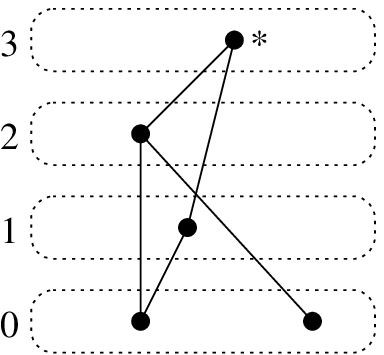}}
\end{minipage}
$$
Next, the $\bob{P}$ element is the single maximal element of the
poset so $x_5=3=\rank(P)-1$ and so we apply rule {\subii}:
$$\mapsto\quad
\begin{minipage}{10em}
  \scalebox{0.70}{\includegraphics{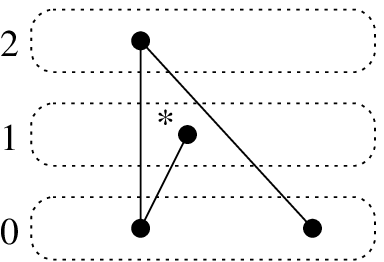}}
\end{minipage}
$$
The marked element is now alone on level $1<\rank(P)-1$ so $x_4=1$
and we apply {\subii} to get:
$$
\mapsto\quad
\begin{minipage}{9.5em}
  \scalebox{0.70}{\includegraphics{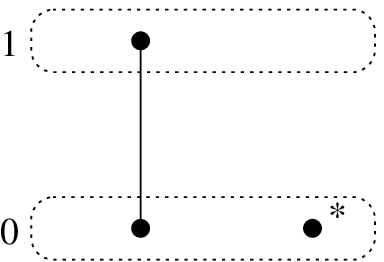}}
\end{minipage}
$$
The marked element is at level 0 so $x_3=0$ and apply {\subi} (since
it is not alone on its level) to get
$$\mapsto\quad
\begin{minipage}{6.3em}
  \scalebox{0.70}{\includegraphics{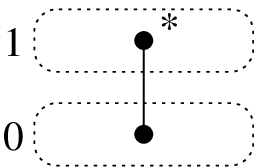}}
\end{minipage} \mapsto\;\;\;
\begin{minipage}{4.5em}
  \scalebox{0.70}{\includegraphics{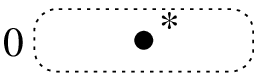}}
\end{minipage}
$$
The final 2 values are easily seen to be $x_2=1$ and $x_1=0$. The
ascent sequence which encodes this {\tpt}-free poset is
$x=(0,1,0,1,3,1,1,2)$.
\end{example}

\section{Main results}\label{mainresults}

For $r \geq 1$, let $G_r(t,u,v,z)$ denote the coefficient of $x^r$
in $G(t,u,v,z,x)$.  Thus $G_r(t,u,v,z)$ is the generating function
of those ascent sequences that begin with $r \geq 1$ 0's followed by
1. We let $G^r_{a,l,m,n}$ denote the number of ascent sequences of
length $n$ which begin with $r$ 0's followed by 1, have $a$ ascents,
 last letter $\ell$, and a total of $m$ zeros. We then let
\begin{equation}
G_r:=G_r(t,u,v,z)=\sum_{a,\ell,m\geq 0,n \geq r+1}
G^r_{a,l,m,n} t^n u^a v^{\ell} z^m.
\end{equation}

Clearly, since the sequence $0\ldots 0$ has no ascents and no
initial run of 0's (by definition), we have that the generating
function for such sequences is $$1+tz +
(tz)^2+\cdots=\frac{1}{1-tz}$$ where 1 corresponds to the empty
word. Thus, we have the following relation between $G$ and $G_r$:
\begin{equation}\label{solution}
G=\frac{1}{1-tz}+\sum_{r\geq 1}G_r \,x^r.
\end{equation}

\begin{lemma}\label{lem2} For $r \geq 1$, the generating function $G_r(t,u,v,z)$ satisfies
\begin{equation}\label{main}
(v-1-tv(1-u))G_r=(v-1)t^{r+1}uvz^r+t((v-1)z-v)G_r(t,u,1,z)+tuv^2G_r(t,uv,1,z).
\end{equation}\end{lemma}

\begin{proof}

Our proof follows the same steps as in Lemma 13
in~\cite{BCDK}. Fix $r \geq 1$.  Let
$x'=(x_1,\ldots,x_{n-1})$ be an ascent sequence
beginning with $r$ 0's followed by 1, with $a$ ascents and $m$ zeros
where $x_{n-1}=\ell$. Then $x=(x_1,
  \dots, x_{n-1}, i)$ is an ascent sequence if and only if $i\in [0,
    a+1]$. Clearly $x$ also begins with $r$ 0's followed by 1. Now,
if $i=0$, the sequence $x$ has $a$
    ascents and $m+1$ zeros. If $1\leq i\leq \ell$, $x$ has
    $a$ ascents and $m$ zeros. Finally if $i \in [\ell+1,a+1]$, then
$x$ has $a+1$ ascents and
    $m$ zeros. Counting the sequence $0\ldots 01$ with $r$ 0's
    separately, we have
\begin{eqnarray*}
G_r &=& t^{r+1}u^1v^1z^r+\sum_{\stackrel{a,\ell,m\geq
0}{n \geq r+1}}G^r_{a,\ell,m,n}t^{n+1}\left(u^av^0z^{m+1}+\sum_{i=1}^{\ell}u^av^iz^m+\sum_{i=\ell+1}^{a+1}u^{a+1}v^iz^m\right)\\
&=& t^{r+1}uvz^r+t\sum_{\stackrel{a,\ell,m\geq
0}{n \geq r+1}}G_{a,\ell,m,n}t^nu^az^m\left(z+\frac{v^{\ell+1}-v}{v-1}+u\frac{v^{a+2}-v^{\ell+1}}{v-1}\right)\\
&=&
t^{r+1}uvz^r+tzG_r(t,u,1,z)+tv\frac{G_r-G_r(t,u,1,z)}{v-1}+tuv\frac{vG_r(t,uv,1,z)-G_r}{v-1}.
\end{eqnarray*}
The result follows.
\end{proof}

Next just like in Subsection 6.2 of ~\cite{BCDK}, we use
the kernel method to proceed. Setting
$(v-1-tv(1-u)) =0$ and solving for $v$, we obtain that
the substitution $v=1/(1+t(u-1))$ will kill the left-hand side of
(\ref{main}). We can then solve for $G_r(t,u,1,z)$ to
obtain  that
\begin{equation}\label{relation}
G_r(t,u,1,z)=\frac{(1-u)t^{r+1}uz^r+uG_r\left(t,\frac{u}{1+t(u-1)},1,z\right)}{(1+zt(u-1))(1+t(u-1))}.
\end{equation}
Next we define
\begin{eqnarray}
\delta_k &=& u-(1-t)^k(u-1) \ \mbox{and} \\
\gamma_k &=& u-(1-zt)(1-t)^{k-1}(u-1)
\end{eqnarray}
for $k \geq 1$. We also set $\delta_0 = \gamma_0 =1$.
Observe that
$\delta_1 = u -(1-t)(u-1) = 1+t(u-1)$ and
$\gamma_1 = u -(1-zt)(u-1) = 1+zt(u-1)$. Thus we can rewrite
(\ref{relation}) as
\begin{equation}\label{relation2}
G_r(t,u,1,z) = \frac{t^{r+1}z^ru(1-u)}{\delta_1 \gamma_1} +
\frac{u}{\delta_1 \gamma_1}G_r(t,\frac{u}{\delta_1},1,z).
\end{equation}

For any function of $f(u)$, we shall write $f(u)|_{u = \frac{u}{\delta_k}}$
for $f(u/\delta_k)$. It is then easy to check that
\begin{enumerate}
\item $ \displaystyle (u-1)|_{u = \frac{u}{\delta_k}} = \frac{(1-t)^k(u-1)}{\delta_k}$,

\item $ \displaystyle \delta_s|_{u = \frac{u}{\delta_k}} = \frac{\delta_{s+k}}{\delta_k}$,

\item $ \displaystyle \gamma_s|_{u = \frac{u}{\delta_k}} = \frac{\gamma_{s+k}}{\delta_k}$, and

\item $ \displaystyle \frac{u}{\delta_s}|_{u = \frac{u}{\delta_k}} = \frac{u}{\delta_{s+k}}$.
\end{enumerate}

Using these relations, one can iterate the recursion (\ref{relation2})
to prove by induction that for all $n \geq 1$,
\begin{eqnarray}
G_r(t,u,1,z) &=& \frac{t^{r+1}z^ru(1-u)}{\delta_1 \gamma_1} +
\left(t^{r+1}z^ru(1-u)\sum_{s=2}^{2^n-1} \frac{u^s(1-t)^s}{\delta_{s}
\delta_{s+1} \prod_{i=1}^{s+1} \gamma_i}\right) + \\
&& \ \ \frac{u^{2^n}}{\delta_{2^n} \prod_{i=1}^{2^n} \gamma_i}G_r(t,\frac{u}{\delta_{2^n}},1,z). \nonumber
\end{eqnarray}

Since $\delta_0 =1$, it  follows that as a power series in $u$,
\begin{equation}\label{GR1}
G_r(t,u,1,z) = t^{r+1}z^ru(1-u)\sum_{s \geq 0} \frac{u^s(1-t)^s}{\delta_s
\delta_{s+1}\prod_{i=1}^{s+1} \gamma_i}.
\end{equation}

We have used Mathematica to compute that
\begin{eqnarray*}
&&G_1(t,u,1,z) = u z t^2+\left(u z+u^2 z+u z^2\right) t^3 \\
&&+ \left(u z+3 u^2 z+u^3 z+u z^2+3 u^2 z^2+u z^3\right) t^4\\
&& +\left(u z+6 u^2 z+7 u^3 z+u^4 z+u z^2+8
u^2 z^2+7 u^3 z^2+u z^3+5 u^2 z^3+u z^4\right) t^5+O[t]^6.
\end{eqnarray*}
For example, the coefficient of $t^4u^2$, $3z+3z^2$ makes sense as
there are 3 ascent sequences of length 4 with 2 ascents and 1 zero,
namely, 0112, 0121, and 0122, while there are 3 ascent sequences of
length 4 with 2 ascents and 2 zeros, namely, 0101, 0102, and 0120
(there are no other ascents sequences of length 4 with 2 ascents).

Note that we can rewrite (\ref{main}) as
\begin{equation}\label{main2}
G_r(t,u,v,z) = \frac{t^{r+1}z^ruv(1-v)}{v\delta_1-1} +
\frac{t(z(v-1)-v)}{v\delta_1-1}G_r(t,u,1,z) +\frac{uv^2t}{v\delta_1-1}G_r(t,uv,1,z).
\end{equation}

For $s \geq 1$, we let
\begin{eqnarray*}
\bar{\delta}_s &=& \delta_s|_{u =uv} = uv-(1-t)^s(uv-1) \ \mbox{and}\\
\bar{\gamma}_s &=& \gamma_s|_{u =uv} = uv-(1-zt)(1-t)^{s-1}(uv-1)
\end{eqnarray*}
and set $\bar{\delta}_0 = \bar{\gamma}_0 = 1$. Then using
(\ref{main2}) and (\ref{GR1}), we have the following theorem.

\begin{theorem}\label{Gr}
For all $r \geq 1$,
\begin{eqnarray}
G_r(t,u,v,z) &=& \frac{t^{r+1}z^ru}{v\delta_1-1}\left(
v(v-1) +t(1-u)(z(v-1)-v)\sum_{s \geq 0} \frac{u^s(1-t)^s}{\delta_s
\delta_{s+1}\prod_{i=1}^{s+1} \gamma_i} \right. \nonumber \\
&& \left. + uv^3t(1-uv)\sum_{s \geq 0} \frac{(uv)^s(1-t)^s}{\bar{\delta}_s
\bar{\delta}_{s+1}\prod_{i=1}^{s+1} \bar{\gamma}_i}\right)
\end{eqnarray}
\end{theorem}

It is easy to see from Theorem \ref{Gr} that
\begin{equation}
G_r(t,u,v,z) = t^{r-1}z^{r-1} G_1(t,u,v,z).
\end{equation}
This is also easy to see combinatorially since every ascent
sequence counted by $G_r(t,u,v,z)$ is of the form $0^{r-1}a$ where
$a$ is an ascent sequence $a$
counted by $G_1(t,u,v,z)$.

We have used Mathematica to compute that
\begin{eqnarray*}
&&G_1(t,u,v,z) = u v z t^2+\left(u v z+u^2 v^2 z+u z^2\right) t^3 \\
&&+\left(u v z+u^2 v z+2 u^2 v^2 z+u^3 v^3 z+u z^2+u^2 z^2+u^2 v z^2+u^2 v^2 z^2+u z^3\right)
t^4\\
&&+\left(u v z+3 u^2 v z+u^3 v z+3 u^2 v^2 z+2 u^3 v^2 z+4 u^3 v^3 z+u^4 v^4 z+u z^2+3 u^2 z^2+u^3 z^2 +3 u^2 v z^2 \right.\\
&&\left.+u^3 v z^2 + 2 u^2 v^2 z^2+2 u^3
v^2 z^2+3 u^3 v^3 z^2+u z^3+3 u^2 z^3+u^2 v z^3+u^2 v^2 z^3+u z^4\right) t^5+
O[t]^6.
\end{eqnarray*}
For example, the coefficient of $t^4u$ is $zv+z^2+z^3$ which makes
sense since the sequences counted by the terms are 0111, 0110, and
0100, respectively.  \\

Note that
\begin{eqnarray*}
G(t,u,v,z,x) &=&  \frac{1}{(1-tz)}+ \sum_{r \geq 1} G_r(t,u,v,z)x^r \\
&=& \frac{1}{(1-tz)}+ \sum_{r \geq 1} t^{r-1} z^{r-1} G_1(t,u,v,z)  x^r \\
&=& \frac{1}{(1-tz)}+ \frac{1}{1-tzx} x G_1(t,u,v,z)
\end{eqnarray*}
Thus we have the following theorem.
\begin{theorem}\label{mainres}
\begin{eqnarray}
&&G(t,u,v,z,x) = \frac{1}{(1-tz)}+ \frac{t^{2}zxu}{(1-tzx)(v\delta_1-1)}
\biggr(  v(v-1) \biggr.  \nonumber \\
&& \biggr.+t(1-u)(z(v-1)-v)\sum_{s \geq 0}
\frac{u^s(1-t)^s}{\delta_s \delta_{s+1}\prod_{i=1}^{s+1} \gamma_i} +
uv^3t(1-uv)\sum_{s \geq 0} \frac{(uv)^s(1-t)^s}{\bar{\delta}_s
\bar{\delta}_{s+1}\prod_{i=1}^{s+1} \bar{\gamma}_i}\biggr).
\end{eqnarray}
\end{theorem}

Again, we have used Mathematica to compute the first few terms of this
series:
\begin{eqnarray*}
&&G(t,u,v,z,x) =  1+z t+\left(u v x z+z^2\right) t^2+\left(u v x z+u^2 v^2 x z+u x z^2+u v x^2 z^2+z^3\right) t^3\\
&&+\left( u v x z+u^2 v x z+2 u^2 v^2 x z+u^3 v^3 x z+u x z^2+u^2 x z^2+u^2 v x z^2\right.\\
&&\left.+u^2 v^2 x z^2+u v x^2 z^2+u^2 v^2 x^2 z^2+u x z^3+u x^2 z^3+u
v x^3 z^3+z^4\right)t^4\\
&&\left(u v x z+3 u^2 v x z+u^3 v x z+3 u^2 v^2 x z+2 u^3 v^2 x z+4 u^3 v^3 x z+u^4 v^4 x z\right.\\
&&+u x z^2+3 u^2 x z^2+u^3 x z^2+3 u^2 v x z^2+u^3 v x z^2+2 u^2 v^2 x z^2+2 u^3 v^2 x z^2+3 u^3 v^3 x z^2\\
&&+u v x^2 z^2+u^2 v x^2 z^2+2 u^2 v^2 x^2 z^2+u^3 v^3 x^2 z^2+u x z^3+3 u^2 x z^3+u^2 v x z^3+u^2 v^2 x z^3\\
&&+u x^2 z^3+u^2 x^2 z^3+u^2 v x^2 z^3+u^2 v^2 x^2 z^3+u v x^3 z^3+u^2 v^2 x^3 z^3+u x z^4\\
&&\left. +u x^2 z^4+u x^3 z^4+u v x^4 z^4+z^5\right) t^5+O[t]^6.
\end{eqnarray*}
One can check that, for instance, the 3 sequences corresponding to
the term $3u^2v^2xzt^5$ are $01112$, $01122$ and $01222$.

\section{Counting \tpt-free posets by size and number of minimal
elements}\label{lenzeros}

In this section, we shall compute the generating function of
\tpt-free posets by size and the number of minimal elements which is
equivalent to finding the generating function for ascent sequences
by length and the number of zeros.

For $n \geq 1$,
let $H_{a,b,\ell, n}$ denote the number of ascent sequences of
length $n$ with $a$ ascents and $b$ zeros which have last letter
$\ell$. Then we first wish to compute
\begin{equation}\label{H1}
H(u,z,v,t) = \sum_{n \geq 1,a,b,\ell \geq 0}
H_{a,b,\ell, n} u^a z^b v^\ell t^n.
\end{equation}

Using the same reasoning as in the previous section, we see
that
\begin{eqnarray*}
H(u,z,v,t) &=& tz+\sum_{\stackrel{a,b,\ell \geq 0}{n \geq 1}}
H_{a,b,\ell,n}t^{n+1} \left(u^av^0z^{b+1}+\sum_{i=1}^{\ell}u^av^iz^b+\sum_{i=\ell+1}^{a+1}u^{a+1}v^iz^b\right)\\
&=& tz+t\sum_{\stackrel{a,b,\ell\geq
0}{n \geq r+1}}H_{a,b\ell,n}t^nu^az^b\left(z+\frac{v^{\ell+1}-v}{v-1}+u\frac{v^{a+2}-v^{\ell+1}}{v-1}\right)\\
&=&
tz+\frac{tv(1-u)}{v-1}H(u,v,z,t)+
\frac{t(z(v-1) -v)}{v-1}H(u,1,z,t) +\frac{tuv^2}{v-1}H(uv,1,z,t).
\end{eqnarray*}
Thus we have the following lemma.
\begin{lemma}\label{lemH1}
\begin{equation}\label{recH1}
(v-1-tv(1-u))H(u,v,z,t) = tz(v-1) +t(z(v-1)-v)H(u,1,z,t) + tuv^2H(uv,1,z,t).
\end{equation}
\end{lemma}
Setting $(v-1-t(1-u)) =0$, we see that the substitution $v =
1+t(u-1) = \delta_1$ kills the left-hand side of (\ref{recH1}). We
can then solve for $H(u,1,z,t)$ to obtain the recursion
\begin{equation}\label{recH2}
H(u,1,z,t) = \frac{zt(1-u)}{\gamma_1} + \frac{u}{\delta_1 \gamma_1}
H(uv,1,z,t).
\end{equation}

By iterating (\ref{recH2}), we can prove by induction that for all
$n \geq 1$,
\begin{equation}\label{recH3}
H(u,1,z,t)= \frac{zt(1-u)}{\gamma_1} + \left( \sum_{s=1}^{2^n-1}
\frac{zt(1-u)u^s(1-t)^s}{\delta_s \prod_{i=1}^{s+1} \gamma_i}\right)
+ \frac{u^{2^n}}{\delta_{2^n} \prod_{i=1}^{2^n} \gamma_i}
H(\frac{u}{\delta_{2^n}},1,z,t).
\end{equation}

Since $\delta_0 =1$, we can rewrite (\ref{recH3}) as
\begin{equation}\label{recH4}
H(u,1,z,t)= \left( \sum_{s=0}^{2^n-1}
\frac{zt(1-u)u^s (1-t)^s}{\delta_s \prod_{i=1}^{s+1} \gamma_i}\right)
+ \frac{u^{2^n}}{\delta_{2^n} \prod_{i=1}^{2^n} \gamma_i}
H(\frac{u}{\delta_{2^n}},1,z,t).
\end{equation}
Thus as a power series in $u$, we can conclude the following.
\begin{theorem}\label{thmH1}
\begin{equation}\label{serH1}
H(u,1,z,t)=\sum_{s=0}^\infty \frac{zt(1-u)u^s (1-t)^s}{\delta_s \prod_{i=1}^{s+1} \gamma_i}.
\end{equation}
\end{theorem}

We would like to set $u=1$ in the power series
$\sum_{s=0}^\infty
\frac{zt(1-u)u^s (1-t)^s}{\delta_s \prod_{i=1}^{s+1} \gamma_i}$,
but the factor $(1-u)$ in the series does not allow us to do that
in this form. Thus our next step is to rewrite the series
in a form where it is obvious that we can set $u=1$ in the series.
To that end, observe that for $k \geq 1$,
\begin{equation*}
\delta_k =u -(1-t)^k(u-1) = 1+u-1 -(1-t)^k(u-1) = 1 -(1-t)^k-1)(u-1)
\end{equation*}
so that

\begin{equation}\label{H4}
\frac{1}{\delta_k} = \sum_{n\geq 0} ((1-t)^k-1)^n(u-1)^n
\sum_{n\geq 0}(u-1)^n \sum_{m=0}^n (-1)^{n-m}\binom{n}{m}(1-t)^{km}.
\end{equation}
Substituting (\ref{H4}) into (\ref{serH1}), we see that
\begin{eqnarray*}
H(u,1,z,t) &=& \frac{zt(1-u)}{\gamma_1} +
\sum_{k \geq 1} \frac{zt(1-u)u^k (1-t)^k}{\prod_{i=1}^{k+1} \gamma_i}
\sum_{n\geq 0} (u-1)^n \sum_{m=0}^n  (-1)^{n-m}\binom{n}{m}(1-t)^{km} \\
&=&  \frac{zt(1-u)}{\gamma_1} + \sum_{n \geq 0} \sum_{m=0}^n
 (-1)^{n-m-1}\binom{n}{m}(u-1)^{n-m}zt \sum_{k \geq 1}
\frac{(u-1)^{m+1}u^k(1-t)^{k(m+1)}}{\prod_{i=1}^{k+1} \gamma_i} \\
&=&  \frac{zt(1-u)}{\gamma_1} + \sum_{n \geq 0} \sum_{m=0}^n
 (-1)^{n-m-1}\binom{n}{m}(u-1)^{n-m}\frac{zt}{(1-zt)^{m+1}} \times \\
&& \ \ \ \sum_{k \geq 1}
\frac{(u-1)^{m+1}(1-zt)^{m+1}u^k(1-t)^{k(m+1)}}{\prod_{i=1}^{k+1} \gamma_i}.
\end{eqnarray*}

Next we need to study the series
$$\sum_{k \geq 1}
\frac{(u-1)^{m+1}(1-zt)^{m+1}u^k(1-t)^{k(m+1)}}{\prod_{i=1}^{k+1} \gamma_i}$$
where $m \geq 0$.
We can rewrite this series in the form
$$-\frac{(u-1)^{m+1}(1-zt)^{m+1}}{\gamma_1} + \sum_{k \geq 0}
\frac{(u-1)^{m+1}(1-zt)^{m+1}u^k(1-t)^{k(m+1)}}{\prod_{i=1}^{k+1} \gamma_i}.$$
We let
\begin{equation}\label{H5}
\psi_{m+1}(u) = \sum_{k \geq 0}
\frac{(u-1)^{m+1}(1-zt)^{m+1}u^k(1-t)^{k(m+1)}}{\prod_{i=1}^{k+1} \gamma_i}.
\end{equation}
We shall show that $\psi_{m+1}(u)$ is in fact a polynomial for all
$m \geq 0$. First, we claim that  $\psi_{m+1}(u)$ salsifies the
following recursion:
\begin{equation}\label{psirec}
\psi_{m+1}(u) = \frac{(u-1)^{m+1}(1-zt)^{m+1}}{\gamma_1} +
\frac{u\delta_1^m}{\gamma_1}\psi_{m+1}\left(\frac{u}{\delta_1}\right).
\end{equation}

That is, one can easily iterate (\ref{psirec}) to prove by induction
that for all $n\geq 1$,
\begin{equation}\label{psirec2}
\psi_{m+1}(u)= \left(\sum_{s=0}^{2^n-1}
\frac{(u-1)^{m+1}(1-zt)^{m+1} u^s(1-t)^{s(m+1)}}{\prod_{i=1}^{s+1} \gamma_i}\right) + \frac{u^{2n} (\delta_{2^n})^m}{\prod_{i=1}^{2^n} \gamma_i}\psi_{m+1}(\frac{u}{\delta_{2^n}}).
\end{equation}
Hence it follows that if $\psi_{m+1}(u)$ satisfies the recursion
(\ref{psirec}), then $\psi_{m+1}(u)$ is given by the power series
in (\ref{H5}).
However, it is routine to check that the polynomial
\begin{equation}\label{psipoly}
\phi_{m+1}(u) = - \sum_{j=0}^m (u-1)^j (1-zt)^j u^{m-j}
\prod_{i=j+1}^m (1 -((1-t)^i)
\end{equation}
satisfies the recursion that
\begin{equation}\label{psipolyrec}
\gamma_1 \phi_{m+1}(u) = (u-1)^{m+1}(1-zt)^{m+1} + u \delta_1^m
\phi_{m+1}\left( \frac{u}{\delta_1}\right).
\end{equation}
Thus we have proved the following lemma.

\begin{lemma}\label{psilem}
\begin{eqnarray}\label{psiform}
\psi_{m+1}(u) &=&  \sum_{k \geq 0}
\frac{(u-1)^{m+1}(1-zt)^{m+1}u^k(1-t)^{k(m+1)}}{\prod_{i=1}^{k+1} \gamma_i}
\nonumber \\
&=& - \sum_{j=0}^m (u-1)^j (1-zt)^j u^{m-j}
\prod_{i=j+1}^m (1 -((1-t)^i).
\end{eqnarray}
\end{lemma}

It thus follows that

\begin{eqnarray*}
H(u,1,z,t) &=&   \frac{zt(1-u)}{\gamma_1} + \sum_{n \geq 0} \sum_{m=0}^n
 (-1)^{n-m-1}\binom{n}{m}(u-1)^{n-m}\frac{zt}{(1-zt)^{m+1}} \times \\
&& -\frac{(u-1)^{m+1}(1-zt)^{m+1}}{\gamma_1} - \sum_{j=0}^m (u-1)^j (1-zt)^j u^{m-j}
\prod_{i=j+1}^m (1 -((1-t)^i).
\end{eqnarray*}
There is no problem in setting $u=1$ in this expression to obtain that
\begin{equation}
H(1,1,z,t) = \sum_{n\geq 0} \frac{zt}{(1-zt)^{n+1}} \prod_{i=1}^n
(1-(1-t)^i).
\end{equation}

Clearly our definitions ensure that
$1+H(1,1,z,t) = P(t,z)$ as defined in the introduction so
that we have the following theorem.
\begin{theorem}\label{mainzeros}
\begin{equation}\label{eqzeros}
P(t,z) = \sum_{n,k \geq 0} p_{n,k}t^n z^k =
1+ \sum_{n \geq 0} \frac{zt}{(1-zt)^{n+1}} \prod_{i=1}^n (1-(1-t)^i).
\end{equation}
\end{theorem}

For example, we have used Mathematica to compute the first few terms
of $P(t,z)$ as
\begin{eqnarray*}
&&P(t,z) = 1+ z t+\left(z+z^2\right) t^2+\left(2 z+2 z^2+z^3\right) t^3+\left(5 z+6 z^2+3 z^3+z^4\right) t^4 \\
&&+\left(15 z+21 z^2+12 z^3+4 z^4+z^5\right) t^5+\left(53 z+84
z^2+54 z^3+20 z^4+5 z^5+z^6\right) t^6+O[t]^7.
\end{eqnarray*}

Next we observe that one can easily derive the ordinary generating
function for the number of \tpt-free posets or, equivalently, for
the number of ascent sequences proved by Bousquet-M\'elou et
al.~\cite{BCDK} from Theorem \ref{mainzeros}. That is, for any
sequence of natural numbers $a=a_1 \ldots a_n$, let $a^+ = (a_1+1)
\ldots (a_n+1)$ be the result of adding one from each element of the
sequence. Moreover, if all the elements of $a=a_1 \ldots a_n$ are
positive, then we let $a^- = (a_1-1) \ldots (a_n-1)$ be the result
of subtracting one to each element of the sequence. It is easy to
see that if $a=a_1 \ldots a_n$ is an ascent sequence, then $0a^+$ is
also an ascent sequence. Vice versa, if $b =0a$ is an ascent
sequence with only one zero where $a = a_1 \ldots a_n$, then $a^-$
is an ascent sequence. It follows that the number of ascent
sequences of length $n$ is equal to the number of ascent sequences
of length $n+1$ which have only one zero.   Hence
\begin{eqnarray*}
P(t) &=& \sum_{n \geq 0} p_nt^n = \frac{1}{t}\frac{\partial P(t,z)}{\partial z}\big|_{z=0} \\
&=& \sum_{n \geq 0} \prod_{i=1}^n (1-(1-t)^i).
\end{eqnarray*}

Results in~\cite{BCDK,cdk,dp} show that \tpt-free posets of size $n$
with $k$ minimal elements are in bijection with the following
objects. (See \cite{BCDK,cdk,dp} for the precise definitions.)
\begin{itemize}
\item ascent sequences of length $n$ with $k$ zeros;
\item permutations of length $n$ {\em avoiding} \pattern whose {\em
leftmost-decreasing run} is of size $k$;
\item {\em regular linearized chord diagrams} on $2n$ points with {\em initial
run of openers} of size $k$; \item {\em upper triangular matrices}
whose non-negative integer entries sum up to $n$, each row and
column contains a non-zero element, and the sum of entries in the
first row is $k$.
\end{itemize}
 Thus (\ref{eqzeros}) provides generating
functions for \pattern-avoiding permutations by the size of the
leftmost-decreasing run, for regular linearized chord diagrams by
the size of the initial run of openers, and for the upper triangular
matrices by the sum of entries in the first row. Moreover,
Theorem~\ref{mainres}, together with bijections
in~\cite{BCDK,cdk,dp} can be used to enumerate the permutations,
diagrams, and matrices subject to 4 statistics. However, we have
chosen not to state explicit generating functions related to the
permutations and diagrams.

Finally, we conjecture that $P(t,z)$ given in
Theorem~\ref{mainzeros} can be written in a different form.

\begin{conjecture}\label{conject}
$$P(t,z) = \sum_{n,k \geq 0} p_{n,k}t^n z^k = \sum_{n \geq 0}\prod_{i=1}^n (1-(1-t)^{i-1}(1-zt)).
$$\end{conjecture}

\section{Restricted ascent sequences and the Catalan numbers}\label{catalan}

Recall\footnote{We would like to thank Anders Claesson for sharing
with us his software to work with posets. Special thanks go to
Hilmar Haukur Gudmundsson for providing us the main ideas, and
essentially a solution to Theorem~\ref{thm2}} that a sequence
$(x_1,\dots , x_n ) \in \mathbb{N}^n$ is a restricted ascent
sequence of length $n$ if it satisfies $x_1=0$ and $x_i \in
[m-1,1+\asc(x_1,\dots , x_{i-1})]$ for all $2\leq i \leq n$, where
$m$ is the maximum element in $(x_1,\dots , x_{i-1})$.

\begin{theorem}\label{thm2} The number of restricted ascent sequences of length $n$ is given by the $n$-th Catalan number.\end{theorem}

\begin{proof} Lets $R_n$ denote the number of restricted
ascents sequences of length $n$.  The Catalan numbers $C_n$
can be defined by the recursion
$$C_{n+1} = \sum_{k=0}^n C_k C_{n-k}$$
with the initial condition that $C_0 =1$. It is easy to see that
$R_0 =1$ since the empty sequence is a restricted ascent sequences.
We must show that
\begin{equation}\label{Rrec}
R_{n+1} = \sum_{k=0}^n R_k R_{n-k}.
\end{equation}
Thus we need a procedure to take a restricted ascent sequence
$D_1$ of length $k$ and a restricted ascent sequence
$D_2$ of length $n-k$ and produce a restricted ascent sequence
$D$ of lenght $n+1$. We shall describe a procedure
``gluing'' two ascent sequences, $D_1$ and $D_2$ which is
equivalent to gluing two Dyck paths together. To
define our gluing procedure we first need the concept of the ``rightmost
maximum'' in an ascent sequence, defined as a left-to-right maximum
$x$ such that $x$ is one more than the number of ascents to the left
of $x$, and none of the left-to-right maxima to the right of $x$ has
this property (in other words, this is the last time we use the
maximum option in the interval $[m-1,1+\asc]$ among the
left-to-right maxima). The sequence $00\ldots0$ is the only one that
does not have the rightmost maximum. For example, 0010101003 has the
rightmost maximum (the leftmost) 1, whereas 0010103323234 has the
rightmost maximum (the leftmost) 3.
Then procedure of ``gluing'' two ascent sequences, $D_1$ and $D_2$ together
can be described as follows.

\begin{enumerate}
\item For $D_1\not=\emptyset$, define $D_1 + D_2 := D_1 (1+asc(D_1))
(D_2$++$(asc(D_1)))$ where ``++'' means increasing {\em each} element of $D_2$ by the number $asc(D_1)$. For example, if
$D_1 = 01021$ and $D_2 =01212$, then $D_1+D_2 = 01021323434$.

\item For $D_1=\emptyset$ define $D_1 + D_2 := D_2$ with the rightmost
maximum element duplicated (add extra 0 if $D_2=00\ldots0$). For
example, $\epsilon+ 01212 = 012212$.
\end{enumerate}

It is easy to see that in Case 1, the element $(1+asc(D_1)$ is the
the rightmost maximum element of $D_1
(1+asc(D_1))(D_2$++$(asc(D_1)))$ which is either the rightmost
element if $D_2 = \epsilon$ or is followed by $asc(D_1)$ if $D_2
\neq \epsilon$ since $D_2$ must start with 0 in that case. It
follows that the rightmost maximal element is not duplicated in
$D_1+D+2$ in Case 1 and, hence, it is easy to recover $D_1$ and
$D_2$ from $D_1 +D_2$.  Clearly, in Case 2, the rightmost maximal
element of $D_1+D_2$ is dulplicated so that we can distinguish Case
1 from Case 2. Moreover, it is easy to see that we can recover $D_2$
from $D_1 +D_2$ in Case 2. This proves that (\ref{Rrec}) holds and
hence $R_n =C_n$ for all~$n$.

Here are examples of decompositions for $n=3$ and $n=4$ ($\epsilon$
stays for the empty word):

\begin{center}
\begin{tabular}{llll}

000 = $\epsilon$ + 00 & 0000 = $\epsilon$ + 000 & 0100 = 0 + 00 & 0112 = 011 + $\epsilon$ \\

001 = 00 + $\epsilon$ & 0001 = 000 + $\epsilon$ & 0101 = 0 + 01 & 0121 = 01 + 0 \\

010 = 0 + 0 & 0010 = 00 + 0 & 0102 = 010 + $\epsilon$ & 0122 = $\epsilon$ + 012 \\

011 = $\epsilon$ + 01 & 0011 = $\epsilon$ + 001 & 0110 = $\epsilon$ + 010 & 0123 = 012 + $\epsilon$\\

012 = 01 + $\epsilon$ & 0012 = 001 + $\epsilon$ & 0111 = $\epsilon$
+ 011 &
\end{tabular}

\end{center}

\end{proof}

 \begin{figure}
        \centering
            \includegraphics[width=100mm]{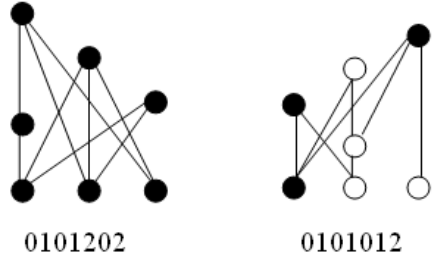}
\vspace{-40mm}
        \caption{Counterexamples to the statement that restricted ascent sequences correspond to
        \tpt- and \tpo-free posets under the bijection in~\cite{BCDK}.}
        \label{fig1}
      \end{figure}

Recall that posets avoiding \tpo\ are those that do not contain, as
an induced subposet, a 3-element chain together with another element
which is incomparable to all elements in the 3-element chain. As we
mentioned in the introduction, the number of posets avoiding \tpt\
and \tpo\ is given by the Catalan numbers
(see~\cite{stanley,SKANDERA}). Using the bijection~\cite{BCDK}
applied to small restricted ascent sequences, one would be tempted
to conjecture that restricted ascent sequences are bijectively
mapped to \tpt- and \tpo-free posets as both of the objects are
counted by the Catalan numbers.  Indeed, this is true
for posets of size less than or equal to six.

Moreover, we can show that the first time one violates that
restricted ascent sequence condition, then the corresponding
\tpt-free poset contains an induced copy of \tpo. That is, suppose
that $a = a_1 \ldots a_n$ is a restricted ascent sequence, $m  =
\max(a_1, \ldots, a_n) \geq 2$, and $x < m-1$. Then we claim that
the poset corresponding to $ax$, must contain an induced copy of
\tpo.  That is, let $r$ be the element on level $x$ that corresponds
to $x$ under the bijection of Section~\ref{decomp}. Now in $ax$, $x$
is preceded by a larger element, and thus $r$ has a neighbor, say
$s$, on its level, level $x$. Because the first time we encounter
$m$ in $a$, its corresponding element $z$ in the poset covers all
maximal elements, it follows that there must be at least one
non-maximal element, say $u$, on level $m-1$. Next, since $x < m-1$,
there exists an element $e$ in the poset such that $e<c$ and $e\not
<b$. That is, $c$ is on a higher level than $b$ and the down-sets
are linearly ordered by inclusion according to their levels. Since
$r$ copies relations of $s$, $e\not<r$. Since $r$ is a maximal
element, also $r\not<e$ and $r \not< u$. Finally, $u$ is a
non-maximal element, thus there exists $v>u$. Finally $v \not< r$
since $r$ is maximal so that the four elements $e<u<v$ and $r$ form
a \tpo\ configuration.

If it was the case that our addition operations preserved the
property of containing \tpo\ configuration, then it would be the
case that the bijection in Section~\ref{decomp} would send \tpt- and
\tpo-free posets to restricted ascent sequences. However, this is is
not the case. For example, consider the poset on the left in
Figure~\ref{fig1} which corresponds to \tpt-free poset corresponding
to the ascent sequence 0101202. One can check that there is an
induced \tpo\ in the poset corresponding to the non-restricted
ascent sequence 010120, but clearly there is no induced \tpo\ in the
\tpt-free poset corresponding to the ascent sequence 0101202. This
means that there must be a restricted ascent sequence of length
seven whose corresponding \tpt-free poset does contain an induced
copy of \tpo.  Such a sequence and its corresponding \tpt-free poset
is shown on the right in Figure~\ref{fig1}.

We leave it as open problem to characterize \tpt-free posets
corresponding to restricted ascent sequences under the bijection
in~\cite{BCDK} and to characterize ascent sequences corresponding to
\tpt- and \tpo-free posets under the same bijection.

\end{document}